\documentstyle[11pt]{article}
\title{On arrangements of real roots of a real polynomial and its derivatives}
\author{Vladimir Petrov Kostov \\ \\ \hspace{7cm}{\sl To Prof. A.A.Bolibrukh}} 
\date{}
\newtheorem{tm}{Theorem}
\newtheorem{defi}[tm]{Definition}
\newtheorem{rem}[tm]{Remark}
\newtheorem{lm}[tm]{Lemma}
\newtheorem{ex}[tm]{Example}

\newtheorem{prop}[tm]{Proposition}
\begin{document}
\maketitle

\begin{abstract}
We prove that all arrangements 
(consistent with the Rolle theorem and some other natural restrictions) 
of the real roots of a real 
polynomial and of its $s$-th derivative are realizable by real polynomials.
\end{abstract}

In the present paper we consider a real polynomial of one real variable 
$P(x,a)=x^n+a_1x^{n-2}+\ldots +a_{n-1}$. We are interested in the question 
what {\em arrangements} between the real roots of $P$ and $P^{(s)}$ are 
possible ($1\leq s\leq n-1$). To define an arrangement means to write down 
the roots of $P$ and $P^{(s)}$ in a chain in which every two consecutive roots 
are connected either by an equality or by an 
inequality $<$. The arrangement $\alpha$ is said to 
belong to the closure of the arrangement $\beta$ if it is obtained from 
$\beta$ by replacing some inequalities by equalitites. The results are the 
first step towards the study of 
real discriminant sets $\{ a\in {\bf R}^{n-1}|{\rm Res}(P,P^{(s)})=0\}$. 

In an earlier paper \cite{KoSh} it is shown that if $P$ is {\em hyperbolic}, 
i.e. with $n$ real roots, then {\em the standard Rolle restrictions} are 
necessary 
and sufficient conditions for a root arrangement to be realizable (see 
Theorems 2 and 4.4 in \cite{KoSh}). Namely, 
denote by $x_1\leq \ldots \leq x_n$ the roots of $P$ and by 
$\xi _1\leq \ldots \leq \xi _{n-s}$ the ones of $P^{(s)}$ (which is also 
hyperbolic). Then one has 

\begin{equation}\label{Rolle}
x_l\leq \xi _l\leq x_{l+s}
\end{equation} 
for $l=1,\ldots ,n-s$ and every arrangement of the roots of $P$ and 
$P^{(s)}$ which is consistent with (\ref{Rolle}) is realizable. One presumes 
also that the following conditions hold:

A) {\em If a root of $P$ of multiplicity $d>s$ coincides with a root of 
$P^{(s)}$ of multiplicity $g$, then $g=d-s$ (self-evident).}
   
B) {\em If a root $\xi$ of $P^{(s)}$ coincides with a root of 
$P$ of multiplicity $\kappa \leq s$, then $\xi$ is a simple root of $P^{(s)}$ 
(see \cite{KoSh}, Lemma~4.2) and one has $\kappa \leq s-1$.}

C) {\em If $x_l=\xi _l$ or $x_{l+s}=\xi _l$, then 
$x_l=x_{l+1}=\ldots =x_{l+s}=\xi _l$ 
(self-evident for $s=1$ and easy to prove by induction on $s$ for $s>1$).}
 
\begin{ex}
If $n=2$, $s=1$, then there 
are two possible arrangements (i.e. consistent with (\ref{Rolle}), A) 
B) and C)) : $x_1<\xi _1<x_2$ and $x_1=\xi _1=x_2$. They are 
both realizable by hyperbolic polynomials. 
\end{ex}

In the present paper we treat the 
case when $P$ is arbitrary (not necessarily hyperbolic). (Notice that 
$P^{(s)}$ can be hyperbolic even if $P$ is not.) 

\begin{defi}\label{defR}
Suppose that $P$ has $m$ conjugate couples of complex 
roots and $n-2m$ real roots. Then a priori $P^{(s)}$ has at least 
$n-2m-s$ real roots counted with the multiplicities. Indeed, a real 
root of $P^{(i)}$ of 
multiplicity $l\geq 1$ is a root of $P^{(i+1)}$ of multiplicity $l-1$ and 
between every two real roots of $P^{(i)}$ there is a root of $P^{(i+1)}$. 
Iterating this rule $s$ times one obtains the existence of  
$n-2m-s$ real roots of $P^{(s)}$ (we call them {\em Rolle} roots) 
which together with the real roots of $P$ satisfy 
conditions (\ref{Rolle}), A) and B). 
A Rolle root is multiple only if it coincides 
with a root of $P$ of multiplicity $>s$. Eventually, 
$P^{(s)}$ can have  
$\leq 2m$ other ({\em non-Rolle}) real roots counted with the multiplicities 
some (or all) of which can coincide with Rolle ones. 
Which real roots of $P^{(s)}$ should be chosen as Rolle and which as 
non-Rolle ones is not always uniquely defined  
and when it is not we assume that a choice is made. 
\end{defi}

\begin{ex}
The polynomial $x^6-x^2=x^2(x^2-1)(x^2+1)$ has real roots $x_1=-1$, 
$x_2=x_3=0$, $x_4=1$ (and complex roots $\pm i$). 
One has $P'=6x^5-2x=2x(\sqrt{3}x^2-1)(\sqrt{3}x^2+1)$, i.e. $P'$ has three 
Rolle roots (and no non-Rolle ones) -- $0$ and $\pm 1/3^{1/4}$ where 
$0$ is a common root 
for $P$ and $P'$, see A). It has also two complex roots $\pm i/3^{1/4}$. 
One has 
$P''=30x^4-2$, i.e. $P''$ has two Rolle roots $\pm 1/15^{1/4}$, no 
non-Rolle ones and two complex roots $\pm i/15^{1/4}$. One has $P'''=120x^3$, 
i. e. $P'''$ has a triple real 
root at $0$ and no complex roots. One copy of this real root should be 
considered as a Rolle one and the other two as non-Rolle ones.
\end{ex} 

\begin{prop}\label{mult}
Suppose that a real root of $P$ of multiplicity $d$ coincides with a 
real root of $P^{(s)}$ of multiplicity $g$. Then 

1) if $d>s$, then one has $g=d-s$; in this case this is a Rolle root of 
$P^{(s)}$ of multiplicity $d-s$; 

2) if $0\leq d\leq s$, then one has $g\leq 2m+1$ (and if $g\geq 1$, 
then $d<s$).
\end{prop}

Observe that in the above example one has $m=1$ and for $s=3$ the estimation 
$2m+1$ is attained by the multiplicity of $0$ as a root of $P'''$. The 
proposition generalizes conditions A) and B) in the case of arbitrary $m$. 

{\em Proof:}

Part 1) is self-evident. Prove part 2). If the root is non-Rolle and does not 
coincide with a Rolle one, then its 
multiplicity is $\leq 2m$. If the root is Rolle and does not coincide 
with a non-Rolle one, then either it coincides with a root of $P$ of 
multiplicity $>s$ and we are in case 1) or it is a simple root. Finally, if 
the root is Rolle and coincides with a non-Rolle one, then the Rolle root must 
be simple (otherwise there will be a contradiction with part 1)) and the sum 
of their two multiplicities is $\leq 2m+1$. ~~~$\Box$

\begin{defi}
An arrangement of the real roots of $P$ and $P^{(s)}$ is called 
{\em a priori admissible} if there exist $n-2m-s$ Rolle roots of $P^{(s)}$ 
in the sense of Definition~\ref{defR} and if conditions 1) and 2) of 
Proposition~\ref{mult} hold.
\end{defi}

\begin{tm}
All a priori admissible root arrangements are realizable by real 
polynomials of degree $n$.
\end{tm}
 
{\em Proof:}

$1^0$. We explain first in $1^0$ -- $7^0$ 
why all a priori admissible arrangements in which the 
derivative $P^{(s)}$ is hyperbolic and which are 
the {\em least generic} are realizable. "Least generic" 
means that all non-Rolle roots of $P^{(s)}$ coincide with 
Rolle ones or with roots of $P$ . The general case is treated 
in $8^0$ -- $11^0$. 

To realize an a priori admissible arrangement with $P^{(s)}$ 
hyperbolic and with the necessary multiplicities of the real roots of $P$ 
consider the family of polynomials 

\begin{equation}\label{P} 
P(x,w,g,t)=\prod _{j=1}^q(x-w_j)^{m_j}\prod _{j=1}^m((x-g_j)^2+t_j^2)
\end{equation} 
where $w_j$, $j=1,\ldots ,q$, are the real roots of $P$, of 
multiplicities $m_j$ ($w_0=0\leq w_1\leq \ldots \leq w_q\leq 1=w_{q+1}$), and 
$g_j\pm it_j$ are its complex 
roots (not necessarily distinct), $t_j\geq 0$, $0\leq g_j\leq 1$. We allow 
here equalities between the roots $w_j$ for convenience; it will be shown 
that the necessary arrangement is realizable for roots with strict inequalities 
between them. 

Denote by $\xi _1\leq \ldots \leq \xi _{n-s}$ the real parts of the 
roots of $P^{(s)}$ ($n-2m-s$ of them are just Rolle roots) and by 
$\theta _1\leq \ldots \leq \theta _m$ the 
biggest nonnegative imaginary parts of the roots of 
$P^{(s)}$ (recall that for a least generic arrangement one has $\theta _j=0$). 
Set $\xi _0=0$, $\xi _{n-s+1}=1$. (Notice that $P^{(s)}$ has not 
more conjugate couples of complex roots than $P$, i.e. not more than $m$.) 
The functions $\xi _i$, $\theta _j$ are 
continuous in $(w,g,t)$. 

$2^0$. Suppose that for the desired arrangement of the real roots of $P$ and 
$P^{(s)}$ the Rolle and non-Rolle roots of $P^{(s)}$ are fixed. Denote the 
non-Rolle roots by $u_1\leq \ldots \leq u_{2m}$. Impose additional 
requirements upon the numbers $g_j$ as follows: if the non-Rolle roots with 
odd indices $u_{2p-1},u_{2p+1},\ldots ,u_{2p+2p'-1}$ belong to the interval 
$[w_j,w_{j+1})$, $j<q$, or to $[w_q,w_{q+1}]$, then we require 
that $w_j\leq g_p\leq \ldots \leq g_{p+p'}\leq w_{j+1}$. 
Define the variables $h_1\leq \ldots \leq h_{q+m}$ as the union of the 
variables $w_j$ ($j=1,\ldots ,q$) and $g_i$ ($i=1,\ldots ,m$) with the order 
defined above. Hence, they belong to the unit simplex $\Sigma _{q+m}$. 

$3^0$. In what follows we assume that the variables $t_j$ belong to some 
interval $[0,N]$ where $N>1$. 
We define with the help of the variables $h_j$, $t_i$ continuous functions 
$\eta _j$, $\zeta _i$ such that 
$(\eta _1,\ldots ,\eta _{q+m})\in \Sigma _{q+m}$, $\zeta _i\in [0,N]$. The 
set ${\cal S}=\Sigma _{q+m}\times [0,N]^m$ is homeomorphic to 
$\Sigma _{q+2m}$. By 
the Brouwer fixed point theorem (see \cite{Do}, p. 57), there exists a 
fixed point of the mapping 
$\tau : {\cal S}\rightarrow {\cal S}$, 
$\tau :(h,t)\mapsto (\eta ,\zeta )$, i.e. a point where one has 
$\eta _j=h_j$, $\zeta _i=t_i$. The functions $\eta _j$, $\zeta _i$ are 
defined such that the arrangement of the real roots of $P$ and $P^{(s)}$ at 
the fixed point is 
the required one.

$4^0$. Define the functions $\eta _j$ by the following rules: 

1) We want to achieve the additional conditions (at the fixed point) 
$g_p=u_{2p-1}$, $\ldots$, $g_{p+p'}=u_{2p+2p'-1}$ 
for all appropriate indices, see $2^0$; therefore we set 
$\eta _{i_1}=\xi _{i_2}$ whenever $h_{i_1}$ is a variable 
$g_{p+l}$ and $\xi _{i_2}$ is the corresponding function 
$u_{2p+2l-1}$;  

2) If a variable $h_j$, which is a root $w_i$ of multiplicity $<s+1$, must 
coincide with a simple root 
$\xi _k$ of $P^{(s)}$ or, more generally, with 
the roots $\xi _k=\xi _{k+1}=\ldots =\xi _{k+l}$, then we set $\eta _j=\xi _k$;

3) If the variables $h_r<h_{r+1}<\ldots <h_{r+l}$
(which are all consecutive roots $w_j$ and 
among which there might be roots $w_j$ of multiplicity $\geq s+1$) lie
between the Rolle roots $\xi _k$ and $\xi _{k+v}$ of $P^{(s)}$ and all roots
among the roots $\xi _{k+1}$, $\ldots$, $\xi _{k+v-1}$ (if $v>1$) coincide 
with roots $w_j$
($r\leq j\leq r+l$) of multiplicity $\geq s+1$, then we set

$\eta _{r+j}=\xi _k+(j+1)(\xi _{k+v}-\xi _k)/(l+2)$, $j=0,1,\ldots ,l$.

\begin{rem}
It follows from rules 1) -- 3) that there are $q+m$ functions $\eta _j$ -- 
as many as the variables $h_j$.
\end{rem}

Recall that the arrangement is least generic, i.e. for every non-Rolle root 
$\xi _i$ of $P^{(s)}$ one has either 
$\xi _i=\xi _{i_1}$ where $\xi _{i_1}$ is a Rolle one or $\xi _i=w_{i_2}=h_j$ 
for some $i_2$, $j$. Denote by $l_1$, $\ldots$, $l_{2m}$ the absolute values  
$|\xi _i-\xi _{i_1}|$ and $|\xi _i-w_{i_2}|$ for all $i$, $i_1$ and $i_2$ as 
above. Set $\Phi =l_1+\ldots +l_{2m}$ and 

\begin{equation}\label{t_i} 
\zeta _i=\left| t_i-\frac{1}{3m}\sum _{j=1}^m\theta _j-
\frac{t_i}{3(N+1)^m}|t_1t_2\ldots t_m-1|-\frac{t_i}{12m}\Phi \right|
\end{equation}

$5^0$. Denote by $t_{i_0}$ the greatest variable $t_i$ at the fixed point 
(see $3^0$). Observe first that one can assume that $t_{i_0}>0$. Indeed, if 
$t_{i_0}=0$, then 
$t_i=0$ for all $i$, $P$ is hyperbolic and the roots of $P$ and $P^{(s)}$ 
define an arrangement $\alpha$ from the closure of the desired 
least generic one $\beta$.

\begin{lm}\label{deformation}
If $t_{i_0}=0$, then there exists a real-analytic 
deformation of $P$ into a real polynomial 
which together with its $s$-th derivative defines the arrangement $\beta$.
\end{lm}

The lemma is proved after the theorem. It allows one to consider only the case 
$t_{i_0}>0$. 

$6^0$. One has  

\[ \zeta _{i_0}=t_{i_0}-\frac{1}{3m}\sum _{j=1}^m\theta _j-
\frac{t_{i_0}}{3(N+1)^m}|t_1t_2\ldots t_m-1|-\frac{t_{i_0}}{12m}\Phi ~.\] 
Indeed, all roots of $P^{(s)}$ lie within the convex hull of all roots of 
$P$ (see \cite{PoSz}, p. 108). Hence, one has $\theta _j\leq t_{i_0}$, 
$j=1,\ldots ,m$. One has 
also $|t_1t_2\ldots t_m-1|\leq t_1t_2\ldots t_m+1< (N+1)^m$ and 
$\Phi \leq 4m$ (because for each term $l_j$ one has $l_j\leq 2$). Thus 

\begin{equation}\label{lessthan} 
\frac{1}{3m}\sum _{j=1}^m\theta _j+
\frac{t_{i_0}}{3(N+1)^m}|t_1t_2\ldots t_m-1|+\frac{t_{i_0}}{12m}\Phi < 
mt_{i_0}/3m+t_{i_0}/3+4mt_{i_0}/12m=t_{i_0}
\end{equation} 
and for $i=i_0$ one can delete the absolute value sign in the right hand-side 
of (\ref{t_i}). But then to have 
$\zeta _{i_0}=t_{i_0}$ one must have $\theta _j=0$ for $j=1,\ldots ,m$, 
$t_1t_2\ldots t_m-1=0$ and $l_1=\ldots =l_{2m}=0$. This means that 
$t_j\neq 0$, i.e. no root $g_j+it_j$ of $P$ will be real, that $P^{(s)}$ will 
indeed be hyperbolic ($\theta _j=0$) and that all non-Rolle roots of 
$P^{(s)}$ equal either roots $w_j$ of $P$ or Rolle roots of $P^{(s)}$.     

\begin{rem}
The condition $N>1$ makes possible the choice of the values of the variables 
$t_i$ so that $t_1t_2\ldots t_m-1=0$. One can prove by analogy with 
(\ref{lessthan}) that $|\zeta _i|<N$, i.e. the mapping $\tau$ is indeed from 
${\cal S}$ into itself.
\end{rem}

$7^0$. A priori the fixed point assures the existence of an arrangement only 
from the closure of the necessary one. The fact that at the fixed point no 
inequality between roots of $P$ is replaced by equality is proved by  
analogy with $4^0$ -- $7^0$ of the proof of Theorem 4.4 from \cite{KoSh} 
where the case of $P$ hyperbolic is considered. The 
proof there shows that equalities replacing inequalities between roots of $P$ 
imply that a root of $P$ of multiplicity $m\geq s+1$ is a root of $P^{(s)}$ 
of multiplicity $\geq m-s+1$ which contradicts part 1) of 
Proposition~\ref{mult}. 
In the general case ($P$ not necessarily hyperbolic) the proof is essentially 
the same, the presence of eventual non-Rolle roots can only increase the 
multiplicity of the root as a root of $P^{(s)}$.
 
Hence, the fixed point provides the necessary arrangement. 

$8^0$. To obtain (in $8^0$ -- $9^0$) 
all arrangements in which $P^{(s)}$ is hyperbolic but which 
are not necessarily least generic we use the same construction but with 
another function $\Phi$. Namely, consider a family of such functions $\Phi$ 
depending on a parameter $b\in ({\bf R}_+,0)$ defined as follows: if instead 
of $\xi _i-\xi _{i_1}=0$, see $4^0$, one must have 
$\xi _i-\xi _{i_1}>0$ or $\xi _i-\xi _{i_1}<0$ (and no root $\xi _j$ or 
$w_j$ lies between $\xi _i$ and $\xi _{i_1}$), then in $\Phi$ we replace the 
absolute value $l_{\nu}=|\xi _i-\xi _{i_1}|$ by $|\xi _i-\xi _{i_1}-b|$ (resp. 
by $|\xi _i-\xi _{i_1}+b|$); in the same way for $\xi _i-w_{i_2}$, see $4^0$. 
In a sense, we obtain the not least generic arrangements by deforming least 
generic ones the deformation parameter being $b$. 

$9^0$. Denote by $F(b)$ the set of fixed points of the mapping 
$\tau$ from $3^0$. 
For $b$ small enough one has $(\eta ,\zeta )\in {\cal S}$. 
The set $F(0)$ contains all limit points of the family of sets $F(b)$ when 
$b\rightarrow 0$ and there exists at least one such limit point because all 
sets $F(b)$ (for $b$ small enough) are non-empty and belong to ${\cal S}$ 
which is compact. Hence, one can choose $b>0$ small enough and a fixed point 
of $F(b)$ at which there is an inequality between two roots in the arrangement 
if there is an inequality in the arrangement for $b=0$, and the equalities 
$\xi _i-\xi _{i_1}=0$ or $\xi _i-w_{i_2}=0$ where this is necessary are 
replaced by the desired inequalities. 

$10^0$. Obtain all arrangements in which $P^{(s)}$ is not hyperbolic and 
which are least generic. Suppose that $P^{(s)}$ must have exactly $m'$ 
conjugate couples of complex roots. In this case we assume that $m'$ of the 
couples of roots $g_j\pm it_j$ are replaced by a couple $\pm iv$ where 
$v>0$ is ``large'', i.e. much bigger than $N$. Hence,  
$P^{(s)}$ also has exactly 
$m'$ couples of conjugate complex roots with ``large'' 
imaginary parts. One has 

\[ Q:=P/v^{2m'}=(1+x^2/v^2)^{m'}\prod _{j=1}^q(x-w_j)^{m_j}\prod _{j=1}^{m-m'}
((x-g_j)^2+t_j^2)~~,\]
i.e. the family $Q$ is a one-parameter deformation of a family of polynomials 
like (\ref{P}) (the role of the small parameter is played by 
$1/v^2$) and the existence of the necessary arrangements can be deduced 
by analogy with $1^0$ -- $7^0$ (see $9^0$ for the role of the small parameter; 
however, the function $\Phi$ is the one from $1^0$ -- $7^0$). 

$11^0$. To obtain the existence of all arrangements (which are 
not necessarily least generic and with $P^{(s)}$ not 
necessarily hyperbolic) one has to 
combine $8^0$, $9^0$ and $10^0$. The theorem is proved.~~~~~$\Box$

{\em Proof of Lemma~\ref{deformation}:}

$1^0$. We assume that $P$ has the same number of distinct real roots as in the 
desired arrangement $\beta$. If not, then one can first deform $P$ within the 
class of hyperbolic polynomials (while remaining in the closure of $\beta$) to 
achieve this condition. See \cite{Ko} for such deformations. 

We begin with two observations: 

1) for $a>0$, $\mu \in {\bf N}\cup \{ 0\}$ and $\nu$ even  
the polynomial $Q=x^{\mu}(x^{\nu}+a)$  
has a $\mu$-fold root for $x=0$ and its $s$-th derivative for $s>\mu$ has a 
$(\mu +\nu -s)$-fold one; $Q$ has also $\nu /2$ couples of conjugate 
complex roots;

2) with $a$, $\mu$ and $\nu$ as above, the polynomial 
$Q_1=x^{\mu}(x^{\nu}+a+aQ_2(x,a))$ where $Q_2$ is a polynomial in $x$ of 
degree $\leq \nu -1$, $Q_2(0,a)\equiv 0$, has $\nu$ complex zeros for $a$ 
small enough and a real 
$\mu$-fold root at $0$; to see this set $a=c^{\nu}$, $x=cy$; one has 
$Q_1(cy,c^{\nu})=c^{\mu +\nu}y^{\mu}(y^{\nu}+1+Q_2(cy,c^{\nu}))$; the 
last polynomial has a $\mu$-fold root at $0$ and $\nu$  
roots which for $c$ 
small enough are close to the roots of $y^{\nu}+1$, hence, are complex.   

$2^0$. Suppose that the polynomial $P$ of degree $n$ 
realizing with $P^{(s)}$ the 
arrangement $\alpha$ has a real root of multiplicity $\mu +\nu$ 
(with $\nu$ even) which 
(in order to obtain the arrangement $\beta$) must split into $\nu /2$ 
couples of conjugate complex roots and into a real root of multiplicity 
$\mu$. (If several roots of $P$ must split, we make them split one by one.) 
Suppose in addition that in the deformed polynomial (denoted by $R$) 
the real root of 
multiplicity $\mu$ must coincide with a root of $R^{(s)}$ of multiplicity 
$\mu +\nu -s$.  
Assume that the bifurcating root is at $0$ and that 

\begin{equation}\label{PP}
P=x^{\mu +\nu}(1+h(x))~~,~~h(0)=0
\end{equation} 
($P$ is not necessarily monic). Construct the necessary 
deformation of $P$ in the form 

\begin{equation}\label{construction} 
R(x,a)=x^{\mu}(x^{\nu}+a+b_{s-\mu}x^{s-\mu}+\ldots +
b_{\nu -1}x^{\nu -1})(1+g(x,a))
\end{equation}
where $a\in ({\bf R},0)$ and $b_i=b_i(a)$ and $g(x,a)$ ($g(0,a)\equiv 0$) 
are defined such that 
all equalities of the form $x_i=\xi _j$ defining the arrangement $\beta$ will 
be preserved. 

$3^0$. Suppose first that in (\ref{construction}) one has $g(x,a)\equiv h(x)$. 
The condition 

\[ (A)~~:~~R^{(s)}~~{\rm has~~a}~~(\mu +\nu -s){\rm -fold~~root~~at~~}0\]
is a 
triangular linear non-homogeneous system with unknown variables $b_i$; the 
system defines unique functions $b_i=b_i^*a$, $b_i^*\in {\bf R}$. This can 
be checked directly.

Suppose that in (\ref{construction}) one has $g=h(x)+\sum _{j=1}^ld_jh_j(x,d)$ 
where $d=(d_1,\ldots ,d_l)\in ({\bf R}^l,0)$ and $h_j$ depend 
smoothly on $d$. Then condition (A) defines 
unique functions $b_i(a,d)=b_i^*a+a\sum _{j=1}^ld_j\tilde{b}_{i,j}(d)$ where 
$b_i^*\in {\bf R}$ and $\tilde{b}_{i,j}$ are smooth in $d$. This can also 
be checked directly.

$4^0$. For each root $w_j\neq 0$ of $P$ of multiplicity $<s$ 
which must be equal to 
a root $\xi _i$ of $P^{(s)}$ denote by $d_j$ the deviation from its 
position in a deformation of $P$. Admitting such deviations means 
that in (\ref{PP}) the function $h$ 
should be replaced by $h(x)+\sum _{j=1}^ld_jh_j(x,d)$. 

Denote by (B) the system of all conditions $w_j=\xi _i$ for all 
such equalities with $w_j\neq 0$ 
characterizing the arrangement $\beta$. 

$5^0$. For 
any deformation  $R^*(x,a,d)=x^{\mu}(x^{\nu}+a+b_{s-\mu}x^{s-\mu}+\ldots +
b_{\nu -1}x^{\nu -1})(1+g(x,d))$ of $P$ (where $b_k$ are considered 
as small parameters) one can find $d$ depending smoothly 
on $a$ and $b_k$ such that for all $a$ small enough all equalities 
from (B) hold. 
This follows from Propositions~11 and 13 from \cite{Ko} where it is shown 
that the linearizations of the conditions (B) w.r.t. $d$ are linearly 
independent. (In \cite{Ko} their linear independence is proved only when 
$P$ is hyperbolic; this independence is an ``open'' property, so it holds for 
all nearby polynomials as well.)

$6^0$. The independence of these linearizations implies that for $a$ small 
enough the system of conditions (B) applied to the deformation 

\[ \tilde{R}(x,a,d)=x^{\mu}(x^{\nu}+a+b_{s-\mu}(a,d)x^{s-\mu}+\ldots +
b_{\nu -1}(a,d)x^{\nu -1})(1+h(x)+\sum _{j=1}^ld_jh_j(x,d))\]
(with $b_i(a,d)$ defined as in $3^0$) defines unique $d_j=d_j(a)$ 
smooth in $a$. Indeed, the linearizations w.r.t. 
$d$ of the system of conditions (B) from $6^0$ and from $5^0$ are the 
same. 

On the other hand, $b_i$ were defined such that condition (A) holds. Hence, 
for $d=d(a)$ and $b_i=b_i(a,d(a))$ (where $a>0$ is small enough) 
the $(\mu +\nu )$-fold root of $P$ at $0$ 
splits into a real $\mu$-fold root at $0$ and $\nu$ complex roots close to 
$0$ (see observation 2) from $1^0$) and $P^{(s)}$ has a 
$(\mu +\nu -s)$-fold root at $0$. The arrangement of the other real roots of 
$P$ and $P^{(s)}$ remains the same.~~~~~$\Box$

{\bf Acknowledgement.} The present text is written under the impetus 
of the fruitful discussions with B.Z. Shapiro, the author's host during a 
short visit at the University of Stockholm. The author 
thanks both him and his university for the kind hospitality.

Author's address: Universit\'e de Nice, Laboratoire de Math\'ematiques, 
Parc Valrose, 06108 Nice Cedex 2, France. e-mail: kostov@math.unice.fr
\end{document}